\documentclass[11pt]{article}
\usepackage{amssymb,psfig,latexsym}

\newtheorem{theorem}{Theorem}[section]
\newtheorem{lemma}{Lemma}[section]
\newtheorem{exam}{Example}[section]

\newcommand{\bsq}{\vrule height .9ex width .8ex depth -.1ex}

\newcommand{\hsp}{\hspace{\parindent}}

\newcommand{\la}{\lambda}

\newcommand{\RR}{{\Bbb R}}

\newcommand{\ZZ}{{\Bbb Z}}
\newcommand{\CC}{{\Bbb C}}
\newcommand{\FF}{{\Bbb F}}

\newcommand{\sC}{{\cal C}}

\newcommand{\sO}{{\cal O}}

\newcommand{\bv}{{\bf v}}
\newcommand{\bw}{{\bf w}}

\newcommand{\bz}{{\bf z}}

\newcommand{\beql}[1]{\begin{equation}\label{#1}}
\newcommand{\eqn}[1]{(\ref{#1})}
\newcommand{\eeq}{\end{equation}}

\catcode`\@=11
\renewcommand{\section}{
        \setcounter{equation}{0}
        \@startsection {section}{1}{\z@}{-3.5ex plus -1ex minus
        -.2ex}{2.3ex plus .2ex}{\large\bf}
        }
\catcode`@=12
\makeatletter
\def\eqalignno#1{\displ@y \ta {\bf s} kip\@centering
  \halign to\displaywidth{\hfil$\@lign\displaystyle{##}$\ta {\bf s} kip\z@skip
    & $\@lign\displaystyle{{}##}$\hfil\ta {\bf s} kip\@centering
    & \llap{$\@lign##$}\ta {\bf s} kip\z@skip\crcr
    #1\crcr}}
\makeatother

\makeatletter
\def\@sect#1#2#3#4#5#6[#7]#8{\ifnum #2>\c@secnumdepth
     \def\@svsec{}\else 
     \refstepcounter{#1}\edef\@svsec{\csname the#1\endcsname.\hskip .75em }\fi
     \@tempskipa #5\relax
      \ifdim \@tempskipa>\z@ 
        \begingroup #6\relax
          \@hangfrom{\hskip #3\relax\@svsec}{\interlinepenalty \@M #8\par}%
        \endgroup
       \csname #1mark\endcsname{#7}\addcontentsline
         {toc}{#1}{\ifnum #2>\c@secnumdepth \else
                      \protect\numberline{\csname the#1\endcsname}\fi
                    #7}\else
        \def\@svsechd{#6\hskip #3\@svsec #8\csname #1mark\endcsname
                      {#7}\addcontentsline
                           {toc}{#1}{\ifnum #2>\c@secnumdepth \else
                             \protect\numberline{\csname the#1\endcsname}\fi
                       #7}}\fi
     \@xsect{#5}}
\def\@theorem#1#2{\it \trivlist \item[\hskip \labelsep{\bf #1\ #2.}]}
\makeatother

\begin{document}
\begin{center}
{\large {\bf Dynamics of  a Family of Piecewise-Linear 
Area-Preserving Plane Maps}} \\
{\large {\bf I. Rational Rotation Numbers}} \\ \bigskip
{\large {\em Jeffrey C. Lagarias}} \\ \smallskip
Department of Mathematics \\
University of Michigan  \\
Ann Arbor, MI 48109-1043 \\
{\tt email:}  {\tt lagarias@umich.edu}\bigskip \\

{\large {\em Eric Rains}} \\ \smallskip
Department of Mathematics \\
University of California-Davis \\
Davis, CA 95616-8633  \\
{\tt email:}  {\tt rains@math.ucdavis.edu}\bigskip \\
\vspace*{2\baselineskip}
(August 2, 2006  version) \\

\vspace*{1.5\baselineskip}
{\bf ABSTRACT}
\end{center}

This paper studies the behavior under iteration of 
the  maps $T_{ab}(x, y) = (F_{ab}(x) - y, x)$ of the
plane $\RR^2$, in which $F_{ab}(x) = ax$ if $ x \ge 0$ and $bx$ if $x < 0$.
The orbits under
iteration correspond to solutions of the nonlinear difference equation
$x_{n+2}= 1/2(a-b)|x_{n+1}| + 1/2(a+b)x_{n+1} - x_n.$
This family of piecewise-linear
maps has the  parameter space $(a, b) \in \RR^2$. 
These maps are area-preserving homeomorphisms of $\RR^2$ that map
rays from the origin into rays from the origin.
The action on rays gives an auxiliary map $S_{ab}: S^1 \to S^1$
of the circle, which
has a well-defined rotation number.
This paper  characterizes the possible dynamics under iteration
of $T_{ab}$ when the auxiliary map
$S_{ab}$ has rational rotation number. It characterizes cases 
where  the map $T_{ab}$ is
a periodic map. 

\vspace*{1.5\baselineskip}

\setlength{\baselineskip}{1.0\baselineskip}

%
%
%
%
%

\section{Introduction}
\hsp
We study the behavior under iteration of the two parameter 
family of piecewise-linear
homeomorphisms of $\RR^2$ given by
\beql{eq101}
T_{ab} (x,y) =
\left\{
\begin{array}{ccc}
(ax-y,x) & \mbox{if} & x \ge 0 , \\
(bx-y,x) & \mbox{if} & x < 0 .
\end{array}
\right.
\eeq
The parameter space is $(a,b) \in \RR^2$.
This map can be written
\beql{eq102}
T_{ab} (x,y) = \left[
\begin{array}{cc}
F_{ab} (x) & -1 \\ 1 & 0
\end{array}
\right]
\left[
\begin{array}{c}
x\\y
\end{array}\right]
 ~,
\eeq
in which
\beql{eq103}
F_{ab} (x) = \left\{
\begin{array}{ccc}
a & \mbox{if} & x \ge 0 , \\
b & \mbox{if} & x < 0;
\end{array}
\right.
\eeq
we view elements of $(x,y)\in\RR^2$ as column vectors.
The formula \eqn{eq102} shows that $T_{ab} (x,y)$ is a homeomorphism, 
since
\begin{eqnarray}\label{eq104}
T_{ab}^{-1} (x,y) & = &
\left[ \begin{array}{cc}
F_{ab} (y) & -1 \\ 1 & 0
\end{array}\right]^{-1} 
\left[\begin{array}{c}x\\y\end{array}\right] 
\nonumber \\
& = & \left[ \begin{array}{cc}
0 & 1 \\ -1 & F_{ab} (y)
\end{array}\right]
\left[\begin{array}{c}x\\y\end{array}\right] 
.
\end{eqnarray}
It preserves the area form $d \omega = dx \wedge dy$, and 
it also maps rays from the origin into rays 
from the origin.

The maps $T_{ab}$ have some claim to being the simplest 
family of nonlinear
area-preserving maps of the plane, which  
gives one reason to study their dynamics under iteration. 
A second reason to consider the dynamics of  $T_{ab}$
(in terms of the parameters $(a, b)$), 
is that iteration of a fixed map  $T_{ab}$  
encodes the solutions of the
second-order nonlinear recurrence
\beql{eq105}
x_{n+2} =  \mu | x_{n+1} |+ \nu x_{n+1}    - x_n ~
\eeq
via
\beql{eq106}
T_{ab} (x_{n+1},x_n ) = (x_{n+2}, x_{n+1} )
\eeq
in which
$$
\mu=   \frac{1}{2}(a - b),~~~~\nu = \frac{1}{2}(a+b). 
$$
This recurrence can be interpreted as a solution to the 
one-dimensional nonlinear difference  equation
of Schr\"{o}dinger type
\beql{eq108}
-x_{n+2} + 2x_{n+1} - x_n + V_\mu (x_{n+1} ) x_{n+1} = E x_{n+1}, ~
\eeq
in which  the ``potential'' $V_{\mu}(x)$ is
\beql{eq109}
V_\mu (x) := \left\{
\begin{array}{ccc}
\mu & \mbox{if} & x \ge 0 , \\
-\mu & \mbox{if} & x < 0.
\end{array}
\right.
\eeq
and the  energy value $E$ is
\beql{eq110}
E := 2- \nu ~.
\eeq
(Strictly speaking, the
potential is $xV_{\mu}(x)$.)
Holding the potential 
$V_\mu$ fixed and letting the parameter $\nu$ vary 
amounts to studying the set of  solutions for
all real energy values  $E = 2 - \nu$.

In this paper and its sequels we  study the orbits of $T_{ab}$ 
for all parameter values,
viewing the results as  giving information about the
set $\Omega_{SB}$ (``semi-bounded'') of parameter 
values having at least one nonzero bounded orbit,
as well as the smaller set $\Omega_{B}$ (``totally bounded'')
of parameter values having 
all orbits bounded. The set $\Omega_{SB}$
can be viewed as an analogue of the fractal ``butterfly'' structure 
of bound state energies appearing in a model of  
Hofstadter \cite{Ho76} for conduction of electrons in
a metal with a stong magnetic field, as explained below.
Our parameters  $(\mu, 2 - \nu)= ( \frac{1}{2}(a-b), 2- \frac{1}{2}(a + b))$
play the role of $(\alpha, E)$ in Hofstadter's model.
In this paper we get information about the set $\Omega_{P}$
of purely periodic maps, and in part II we obtain information
about the set $\Omega_{R}$ of maps topologically conjugate to a rotation of
the plane, in which case the map has invariant circles. We have
$$
\Omega_{P} \subset \Omega_{R} \subseteq \Omega_{B} 
\subseteq \Omega_{SB},
$$
it is known that $\Omega_{R} = \Omega_{B}$, as discussed
below.
We defer a detailed study of $\Omega_{SB}$ to part III,
where we  show that $\Omega_{SB}$ is a closed set 
which has certain properties like those ascribed to the 
Hofstadter ``butterfly''. 

The dynamics of $T_{ab}$ under iteration includes the induced
dynamics on rays, which is described by a 
circle map $S_{ab}$ (defined in \S2),
as well as a stretching/shrinking  motion within rays.
As is usual for circle maps, 
there is a  dichotomy in the dynamics
between the cases of rational rotation
number and irrational rotation number.
In this paper the main
result is to  classify  the possible dynamics in case of  rational
rotation number. In the process we  characterize  
parameter values where the motion is purely
periodic; this  was done earlier 
by Beardon, Bullett and Rippon \cite{BBR95},
as discussed below. 
Our interest in  purely periodic
maps $T_{ab}$  was stimulated by the observation of
Morton Brown \cite{Br83} in 1983 that
the recurrence 
$x_{n+1} = |x_n| - x_{n-1}$ has
all solutions periodic of period $9$; 
this is the case $\mu =1$, $\nu = 0$ of  \eqn{eq105}.
At the end of the paper we  establish a few facts about the
irrational rotation number case, which is considered
in more detail in part II. In \S2 we summarize the main
results, and then establish them in \S3-\S5.

There has been  substantial earlier work on 
the dynamics of these maps. In 1986 M. Herman
\cite[Chap. VIII]{He86} studied 
the dynamics of maps in the {\it Froeschl\'{e} 
group} $\FF$, which he 
defined to be (\cite[p. 217]{He86})
the  group of homeomorphisms of the plane
generated by the elements of $SL(2, \RR)$, together with
the piecewise linear maps of the form
$$
G_{ab} = \left[
\begin{array}{cc}
1 & 0 \\ F_{ab}(x) & 1
\end{array}
\right]
$$
for real $a,b$, viewed as acting on column vectors.
This group was named after Froeschl\'{e} \cite{Fr68}, who
numerically computed (apparant) invariant circles for
certain maps in $\FF$.
The maps $T_{ab}$ belong to $\FF$, since
$$
\left[ \begin{array}{cc}
F_{ab}(x) & -1 \\ 1 & 0
\end{array}\right] 
= 
\left[ \begin{array}{cc}
-F_{-a, -b}(x) & -1 \\ 1 & 0
\end{array}\right] 
= 
\left[ \begin{array}{cc}
0 & -1\\ 1& 0
\end{array} \right]
 \left[ \begin{array}{cc}
1 & 0 \\ F_{-a, -b}(x)& 1 \end{array}
\right].
$$
Herman was in large part concerned with proving
the existence of invariant circles for a large set of
such maps, and we comment further on his work in part II.
However he also 
conjectured the existence of maps in the Froeschl\'{e} group having
dense orbits in the plane (\cite[p. 221]{He86}).
It remains an open question whether any maps $T_{ab}$ have
this property.
 
In 1995   Beardon, Bullett and Rippon \cite{BBR95}  studied
periodic orbits of the maps \eqn{eq105}.
They also studied the set of parameter values $\Omega_{B}$ 
for which all orbits are bounded, determined restrictions
on the parameters in this set, and observed that the equality
$\Omega_{R}=\Omega_{B}$  follows from results of 
Herman \cite[VIII.2.4]{He86}.
Our methods and results in part I 
overlap considerably with theirs;
more remarks on this are made in
\S2. The results and examples  given here  on periodic orbits
were first obtained  in their paper. We include our own 
 proofs of these results for the reader's
convenience, as our notation differs from that
of \cite{BBR95}, and because properties of certain
families of periodic maps are needed in later proofs
in parts II and III. Our main result here settles 
a question raised in their paper.

We conclude this introduction with two
remarks. First, we note that 
another extensively studied two-parameter family of
plane maps which are piecewise-affine with two pieces are 
the Lozi maps,
introduced by Lozi~\cite{Lo79}.
These maps have the form 
$$L(x, y) = ( 1- a|x| + y, bx),$$ 
where $a, b$ are real and $b <0$. These maps
are  homeomorphisms of $\RR^2$, but do not preserve area 
(except for $b=-1$), nor do they  take rays through the origin to rays
through the origin.  Lozi maps  are known to exhibit
a wide range of chaotic behavior, including strange
attractors, see Misiurewicz~\cite{M79} and Collet and Levy~\cite{CL84}. 
The maps we consider in this paper are much simpler than
the Lozi maps, but our results  show
they have nontrivial behavior under iteration.  

Second, we note that the difference equation \eqn{eq108} may be compared 
with two much-studied classes of operators.

(i) The nonlinear Schr\"{o}dinger operator with a general 
(real-valued) potential
$V$ takes the form
$$
-\Delta u + V (|u| ) u = i\frac{\partial}{\partial t}u, ~
$$
where the function  $u(t, x)$ is
complex-valued. For stationary phase
solutions $u(t, x)= e^{iEt}z(x)$ the phase component factors out
and leads to 
$$
-\Delta z + V (|z| ) z = -Ez. 
$$
The discretized version of this equation involves the
{\em discrete nonlinear Schr\"{o}dinger operator}
$$
\Phi_V(z)_n = -z_{n+2} + 2z_{n+1} - z_n - V (|z_{n+1}| ) z_{n+1} 
$$
in the eigenvalue equation
 \beql{108c}
-z_{n+2} + 2z_{n+1} - z_n - V (|z_{n+1}| ) z_{n+1} = E z_{n+1}, ~
\eeq
for $z_n \in \CC$, see \cite{WS90}.
 The special case of solutions to this operator 
with constant phase $\bz_n = e^{\i \phi} x_n$ (with 
$0 \le \phi < \pi,$ $x_n$ real),
reduces to the general form 
\eqn{eq108} 
 with a symmetric potential  $V(x)= V(-x)$
(by cancelling out all phases).
The recurrence \eqn{eq108} that we consider has
a similar form to \eqn{108c}, with the difference that the potential 
$V(x)$ is  antisymmetric, i.e. 
$V(x) = -V(-x)$.

(ii) {\em  Discrete linear Schr\"{o}dinger operators} on the line have 
been extensively studied,  as simple models for 
conduction/insulation transitions
of elections in metals, and more recently in  quasicrystals,
see \cite{Be92}, \cite{DP86}, \cite{Ho76},
\cite{Ko90},
\cite{Sp86}, \cite{SK87}, \cite{Su95}.
The linear difference operator 
\beql{eq111b}
-x_{n+2} + 2x_{n+1} - x_n + V (n+1 ) x_{n+1} = E x_{n+1}, ~
\eeq
with a potential $V(n)$ depending only on position
is often called the ``tight-binding'' approximation
to the Schr\"{o}dinger operator on the line.
Here one is interested in characterizing 
the values of $E$ which have {\em extended states}
these are bounded orbits, that is,  real-valued orbits in $l_{\infty}(\ZZ)$.
(Extended states are required in the physics literature to not 
belong to  $l_2$ but we do not
impose  this condition; states belonging to $l_2$ are called
 localized states.)
We define the  {\em $l_{\infty}$-spectrum} 
for a  fixed potential $V$  to be 
$$
Spec_{\infty}[V] := \{ E:~ \mbox{The discrete Schr\"{o}dinger equation}~
\eqn{eq111b}~~ \mbox{has a bounded orbit} \}.
$$
In 1976 Hofstadter \cite{Ho76} computed  numerically 
the $l_{\infty}$-spectrum $\Sigma_{\lambda, \alpha}$ 
for the discrete
Schr\"{o}dinger operator with a quasiperiodic potential
$V(n) = \lambda \cos (2 \pi \alpha n )$,  
holding  $\lambda=2$ fixed, and letting  $\alpha$ vary.
(This equation models allowable conduction
energies of an  electron moving in a two-dimensional cubic crystal
with a strong magnetic field applied perpendicularly.)   
He observed that the spectrum formed a two dimensional picture
in the $(\alpha, E)$ plane 
resembling a fractal ``butterfly''. He gave a
a conjectural  explanation for the fractal structure using a kind of
renormalization. For fixed irrational $\alpha$
the $l_{\infty}$-spectrum appeared to be a Cantor set of measure zero.
Hofstadter's model has been extensively studied, with one
goal being to  give a rigorous justification 
of this ``butterfly'' structure; see Bellissard \cite{Be94}
and Sj\"{o}strand \cite{Sj91}.
This spectrum for fixed parameters $(\lambda, \alpha)$
has a ``band'' structure when $\alpha$ is rational. 
A slightly more general problem 
almost Mathieu equation, which 
considers the potential $V(n) =  2\lambda \cos(2\pi((n\alpha + \theta))$
where $(\lambda, \alpha, \theta)$ are parameters.
Much recent progress has been made on the $l^2$-spectrum  of this
equation as a function of the parameters, see Jitomirskaya \cite{Ji99}
and Puig \cite{Pu04}. A solution to  
long-standing problem of Cantor
set spectrum when $\alpha$ is irrational has apparantly been achieved
by  Avila and Jitomirskaya \cite{AJ05}.
For general results on discrete
Schr\"{o}dinger operators see Bougerol and Lacroix\cite{BL85}
and Pastur and Figotin \cite{PF92}.

\noindent \paragraph{Notation.} We write $\bv=(\bv_x, \bv_y) \in \RR^2$,
to be viewed as a column vector.
 An interval $[\bv_1, \bv_2)$ of the 
unit circle, or corresponding sector $\RR^{+}[\bv_1, \bv_2)$
of the plane $\RR^2$, is the one specified by going counterclockwise from
$\bv_1$ to $\bv_2$.

\noindent \paragraph{Acknowledgments.} We did most of
the work on this paper while employed at AT\&T Labs-Research,
whom we thank for support; most of the results in parts I and II
were obtained in the summer of 1993. We thank T. Spencer for
helpful comments on the relation of 
\eqn{eq108} to nonlinear Schr\"{o}dinger
operators, and  M. Kontsevich for bringing 
the work of Bedford,
Bullett and Rippon \cite{BBR95}
to our attention.
%
%
%
%
%

\section{Summary of Results}
\hsp
The parameter space of the map can be taken to be either
$(a, b)$ or $(\mu,\nu)$, as these are equivalent by 
\beql{eq201}
\mu = \frac{1}{2}( a - b), \qquad \nu= \frac{1}{2} (a + b).
\eeq
Both coordinate systems have their advantages, and we
 write the map \eqn{eq101} as $T_{ab}$, $T_{\mu \nu}$ accordingly.
It is convenient to represent the action of $T_{ab}$, 
acting on column  vectors $\bv_n = (x_{n+1}, x_n )$ as
\beql{eq202}
T_n (\bv_0) = \left[ {x_{n+1} \atop  {x_n}} \right]
 = M_n (\bv_0) 
\left[ {x_1 \atop  x_0} \right] ,
\eeq
in which
\beql{eq203}
M_n (\bv_0) = \prod_{i=1}^n \left[
\begin{array}{cc}
F_{ab} (x_i) & -1 \\ 1 & 0
\end{array}
\right] :=
\left[
\begin{array}{cc}
F_{ab} (x_n) & -1 \\ 1 & 0
\end{array} 
\right]
\cdots
\left[
\begin{array}{cc}
F_{ab} (x_2 ) & -1 \\ 1 & 0 
\end{array}
\right] 
\left[
\begin{array}{cc}
F_{ab} (x_1 ) & -1 \\
1 & 0
\end{array}
\right]\,.
\eeq

Conjugation by the involution $J_0: (x,y) \to (-x, -y)$ gives
\beql{eq204}
T_{ba} (x,y) = J_0^{-1} \circ T_{ab} \circ J_0 ~.
\eeq
Thus, in studying dynamics, without loss of generality we can
restrict to the closed half-space
$\{(a,b) : a \ge b \}~$ of the $(a,b)$ parameter space.
This corresponds to the region $\{(\mu, \nu): \mu \ge 0 \}$
of the $(\mu. \nu)$ parameter space, with $T_{\mu \nu}$ conjugate
to $T_{-\mu, \nu}$.

Conjugation by the involution $R: (x,y) \to (y, x)$ takes
$T_{ab}$ to its inverse map $T_{ab}^{-1}$, i.e.
\beql{eq204a}
T_{ab}^{-1} (x,y) = R^{-1} \circ T_{ab} \circ R ~,
\eeq
see Theorem~\ref{th34}.
Thus, the dynamics of running the iteration backwards is 
essentially the
same as running it forwards. 

The map $T_{ab}$ is homogeneous, so  sends rays
$[\bv] := \{ \lambda \bv : \la \ge 0 \}$ to rays $T_{ab} ([\bv])$.
Therefore we obtain a well-defined circle map $S_{ab} : S^1 \to S^1 ~,$
for $0 \le \theta \le 2 \pi,$ given by
\beql{eq206}
S_{ab} (e^{i \theta} ) :=
\frac{T_{ab} (e^{i\theta})}{\|T_{ab} (e^{i \theta}) \|} ~,\qquad
0 \le \theta \le 2 \pi \,.
\eeq
where $\| T_{ab} (e^{i \theta} ) \|$ is the Euclidean norm on 
$\RR^2$, and we identity $e^{i \theta} = x+ yi \in \CC$ with
$(x,y) \in \RR^2$.
By abuse of language we shall also sometimes treat the
circle map as $S_{ab}(\theta)$ having domain $\RR/2\pi\ZZ.$
Understanding the  dynamics of the map $T_{ab}$ subdivides 
into two problems: 
study of the dynamics of the circle map $S_{ab}$, 
and study of the motion of points inside the individual rays.

In \S3 we study the circle map $S_{ab}$.
It  has a well-defined 
rotation number $r(S_{ab} )$, and 
we determine the allowed range of the rotation number, as follows.
%
%
\begin{theorem}\label{th21}
 For fixed real $a$, and $-\infty< b < \infty$,
the rotation number $r(S_{ab})$ is continuous and
nonincreasing in $b$, and completely fills out the
following intervals.

(i) For $a < 0$,
\beql{eq207}
r(S_{ab} ) \in \left[ 0, \frac{1}{2} \right] ~.
\eeq

(ii) Let $0 \le a < 2$. Then for each integer $n \ge 2$,
on the interval $2\cos \frac{\pi}{n} \le a < 2\cos \frac{\pi}{n+1},$
\beql{eq208}
r(S_{ab} ) \in \left[ 0, \frac{1}{n+1} \right].
\eeq

(iii) For $a \ge 2$ one has
\beql{eq209}
r(S_{ab} ) =0 ~.
\eeq
\end{theorem}

We next determine the range of the rotation number in
$(\mu, \nu)$-space, for fixed $\mu$.
%
%
\begin{theorem}\label{th22}
For fixed real $\mu$ and $- \infty \le \nu < \infty$ 
the values of $r(S_{\mu \nu} )$ are continuous and
nonincreasing in $\nu$ and completely fill out the
closed interval
\beql{eq209b}
r(S_{\mu \nu} ) \in \left[ 0, \frac{1}{2} \right] ~.
\eeq
\end{theorem}

We conclude \S3 with examples of specific parameter sets where
the rotation number can be determined exactly. These
include one-parameter families, and cases  where
the map $T_{ab}$ is purely periodic. These examples were
all given  in Bedford, Bullett and Rippon \cite{BBR95}. 

In \S4 we study maps $T_{ab}$ for which $r(S_{ab} )$ is rational.
The maps $S_{ab}$ exhibit a ``mode-locking''
behavior so that there are open sets in the parameter space
for which $r(S_{ab} )$ takes a fixed rational value
(for certain rationals), as in Theorem~\ref{th21}(iii) above.
We first characterize those values of $(a,b)$ where $T_{ab}$ is periodic.

%
%
\begin{theorem}\label{th23}
$T_{ab}$  is of finite order if and only if the orbit of $(0,1)$ 
is periodic.
\end{theorem}

This result was obtained in Beardon, Bullett and Rippon 
\cite[Theorem 3.1(i)]{BBR95}.
This result implies that the parameter values at which
$T_{ab}$ is a periodic map fall in a countable number of
one-parameter families, plus a countable number of isolated
values, described by the period and possible symbolic dynamics
of the orbits, as explained in \S4.

We then prove the  main result of this paper,
which characterizes the dynamics in all cases of rational
rotation number, as follows.
%
%
\begin{theorem}\label{th24}
If the rotation number
$r(S_{ab} )$ is rational, then $S_{ab}$ has a periodic orbit,
and one of the following three possibilities occurs.

(i) $S_{ab}$ has exactly one periodic orbit.
Then $T_{ab}$ has exactly one periodic orbit (up to scaling) 
and all other orbits diverge in modulus to $+ \infty$ as $n \to \pm \infty$.

(ii) $S_{ab}$ 
has exactly two periodic orbits.
Then $T_{ab}$ has no periodic orbits.
All orbits of $T_{ab}$ diverge in modulus
to $+\infty$ as $n \to \pm \infty$, with  
the exception of orbits lying over the two periodic orbits of $S_{ab}$.
These exceptional orbits have  modulus diverging  to $+ \infty$ 
in one direction and to $0$ in the other direction, with forward
divergence for one, and backward divergence for the other.

(iii) $S_{ab}$ has at least three periodic orbits.
Then $T_{ab}$ is of finite order, i.e. $T_{ab}^{(k)} =I$ for some $k \ge 1$,
and all its orbits are periodic.
\end{theorem}

This result 
answers a question raised in Beardon, Bullett and
Rippon \cite[p. 671]{BBR95}: Can there be a nonperiodic $T_{ab}$
having two disjoint orbits of rays on each of which $T_{ab}^q= I$?
This corresponds to case (ii) above, so the answer is: no.

In \S5 we begin the study of maps $T_{ab}$ for which the rotation number 
 $r(S_{ab} )$ is irrational.
We show that the existence of an orbit with elements of modulus bounded 
away from 0 and $\infty$ implies that there exists an invariant circle
(Theorem~\ref{th52}).
Then we prove that an invariant circle is 
necessarily preserved under the
reflection symmetry $R(x,y)=(y,x)$ (Theorem~\ref{th52a}).
Part II makes a further study of maps with irrational rotation number
and invariant circles. We note that
in the irrational rotation number case
 a result of Herman \cite[VIII.2.4]{He86} implies a dichotomy:
either $T_{ab}$ is topologically conjugate to a rotation of
the plane, or else it has a dense orbit. However in part III
we show  in the latter case that it also possesses a bounded orbit.

As remarked in the introduction,
our methods and results overlap with those
used in Beardon, Bullett and Rippon \cite{BBR95}.
Their paper contains additional results about periodic
maps beyond those given here. 
To compare results, their parameters
$(\lambda, \mu)$ correspong to our $(\mu, \nu)$; they lift the 
difference equation  \eqn{eq105} to a map of the plane
exchanging $x$ and $y$ coordinates from
our map, so their rotation numbers are negative
while ours are positive; their parameter sets 
$\Lambda_{P}, \Lambda_{R}, \Lambda_{P}$ correspond to our parameter sets
$\Omega_{P}, \Omega_{R}, \Omega_{B}$, respectively.
Theorem~\ref{th24} is able to advance beyond 
their results using the uniqueness assertion
in Lemma~\ref{le41}.

%
%
%
%
%

\section{Associated Circle Map}
\hsp
We study the circle map 
$S_{\mu \nu}(\theta)$ in the $(\mu, \nu )$-coordinates.
%
%
\begin{theorem}
\label{th31}
Each $S_{\mu\nu}(\theta)$ is an orientation-preserving homeomorphism
from $S^1$ to itself.
The derivative $\frac{d}{d\theta} (S_{\mu\nu} (\theta ))$ 
is continuous and of bounded variation.
\end{theorem}

\paragraph{Proof.}
Coordinatizing $S^1$ as $(\cos(\theta),\sin(\theta))$,
we have
$$
S_{\mu\nu}(\theta)=
\cot^{-1}({\nu\cos(\theta)+\mu|\cos(\theta)|-\sin(\theta)\over \cos(\theta)})+
\pi\left[\cos(\theta)<0\right].
$$
This is clearly continuous, mod $2\pi$; the only cause for concern is when
$\cos(\theta)=0$, where it can be easily verified that the pieces match up.
Similarly, if we look at the derivative of $S_{\mu\nu}$, we get:
$$
{d\over d\theta}S_{\mu\nu}(\theta)={1\over
(\nu\pm\mu)^2\cos^2(\theta)-2(\nu\pm\mu)\cos(\theta)\sin(\theta)+1},
$$
with $+$ when $\cos(\theta)>0$ and $-$ when $\cos(\theta)<0$.
This is piecewise continuous, and again, the pieces match up.  Finally,
the second derivative of $S_{\mu\nu}$ exists except when $\cos(\theta)=0$, 
and is bounded. We have
$$
{d^2\over d\theta^2}S_{\mu\nu}(\theta)=\frac{(\nu\pm\mu)^2\sin(2\theta)+
2(\nu\pm\mu)\cos(2\theta)}
{((\nu\pm\mu)^2\cos^2(\theta)
-2(\nu\pm\mu)\cos(\theta)\sin(\theta) +1)^{2}},
$$
and the numerator is bounded by 
$$
|(\nu\pm\mu)^2\sin(2\theta)+2(\nu\pm\mu)\cos(2\theta)|\le 
(|\nu|+|\mu|)\sqrt{(|\nu|+|\mu|)^2+4}$$
and the denominator by
$$
|(\nu\pm\mu)^2\cos^2(\theta)-2(\nu\pm\mu)\cos(\theta)\sin(\theta)+1|^{-1}
\le 1+{1\over 2}((|\nu|+|\mu|)^2+(|\nu|+|\mu|)\sqrt{(|\nu|+|\mu|)^2+4}).
$$
Thus, we can conclude that ${d\over d\theta}S_{\mu\nu}(\theta)$ is continuous
and of bounded variation.  Finally, $S_{\mu\nu}$ is orientation-preserving,
since ${d\over d\theta}S_{\mu\nu} (\theta )$ is always positive, 
and is a homeomorphism,
since $S^{-1}_{\mu\nu}$ obtained by conjugating $S_{\mu\nu}$ by reflection
through the line $x=y$.~~~$\bsq$ \\

The dynamics of iteration of maps on the circle is well understood.
A basic result (\cite[p.33]{Ni71}) is that every
orientation-preserving homeomorphism $S: S^1\to S^1$ has a well-defined
{\em rotation number}, defined by
\beql{eq201A}
r(S) :=\lim_{n\to\infty} {\tilde{S}^{(n)}(x)-x\over n},
\eeq
where $\tilde{S}$ is any lift of $S$ to ${\bf R}$ ($2\pi$ lifts to $1$);
here $r(S)$ is
independent of the initial value $x$ $(\bmod~1)$ and of the choice of lift.
(Positive rotation number corresponds to counterclockwise rotation.)  
In particular, if
$S$ has a periodic point, then $r(S)$ is rational (and conversely, as well).

This section gives general results
on the rotation number $r(S_{ab} )$; the following
two sections cover the case $r(S_{ab})$ rational and $r(S_{ab})$ irrational,
respectively.
\begin{theorem}\label{th32}
{\rm (i)} For fixed $b$, the  rotation number $r(S_{ab} )$ 
is nonincreasing in $a$, and for fixed $a$ 
it is nonincreasing in $b$.

{\rm (ii)}
For fixed $\mu$, the rotation number $r(S_{\mu \nu})$ 
is nonincreasing in $\nu$.
\end{theorem}

\paragraph{Proof.}
(i). It suffices to show that $r(S_{ab} )$ is 
nonincreasing in $a$,
the result for $b$ follows from  the relation $r(S_{ab})= r(S_{ba})$.

Now, consider the behavior of $S_{ab}(\theta)$ as $a$ increases.
If we can show that it moves clockwise as $a$ increases, for every $\theta$,
then we are done, by inspection of equation \eqn{eq201A}.
To show that $S_{ab}(\theta)$ moves clockwise with $a$, 
we need only show that
$$
L=
T_{ab}(\bv)_x{\partial\over\partial a}T_{ab}(\bv)_y -
T_{ab}(\bv)_y{\partial\over\partial a}T_{ab}(\bv)_x \le 0,
$$
for all $\bv\in \RR^2$.  For $\bv_x<0$, this is clearly $0$,
since $T_{ab}(\bv)$ is independent of $a$ in that case.  
For $\bv_x\ge 0$, $T_{ab}(\bv)=(a \bv_x- \bv_y, \bv_y)$; 
a simple calculation gives
$$
L=-\bv_x^2\le 0,
$$
so (i) is proved.

(ii). As $\nu$ increases, both $a$ and $b$ increase, 
so the result  essentially follows from (i). More precisely, for
fixed $\theta$ the point
$S_{\mu\nu}(\theta)$ moves clockwise as $\nu$ increases.~~~$\bsq$ \\

It a natural quesion to ask how the rotation number 
varies when $\nu$ is held fixed, and $\mu$ varies.
Numerical  evidence indicates that it
is nonincreasing for $\mu >0$, hence nondecreasing for $\mu < 0$
(recall that $r(S_{\mu \nu}) = r(S_{-\mu, \nu})$
since $T_{\mu \nu}$ is topologically conjugate to $T_{-\mu, \nu}$).
This was observed by Bedford, Bullett and Rippon \cite{BBR95},
who proved this holds for $\nu=0$.
One can show  that the 
rotation number is nonincreasing
for $\mu > |\nu|$, using an extension of the  approach used  in
Theorem~\ref{th32}(i). One checks 
for fixed $\theta$ that the fourth iterate
$S_{\mu\nu}^{(4)}(\theta)$ moves clockwise as $\mu$ increases,
with $\nu$ held fixed; there are 16 possible cases to consider,
depending on the signs of the iterates. However 
new ideas seemed needed to determine
how the rotation number behaves in the region  $0 < \mu < |\nu|$.

The fact that the rotation number is monotonic in $b$ allows us to
give fairly strong bounds on the rotation number, for a given $a$,
as follows.

%
%

\begin{theorem}\label{th33}
For each fixed $a$ and $- \infty \le b \le \infty$ 
the values of $r(S_{ab} )$ completely fill out the following intervals,
including  their endpoints.

(i) For $a < 0$, the rotation number satisfies
\beql{eq305}
r(S_{ab} ) \in \left[0, \frac{1}{2}\right] ~.
\eeq

(ii) Let $0 \le a \le 2$. For each $n \ge 2$ on the interval
$2\cos \frac{\pi}{n} \le a < 2\cos \frac{\pi}{n+1}$,
the rotation number satisfies 
\beql{eq306}
r(S_{ab} ) \in \left[0, \frac{1}{n+1} \right].
\eeq

(iii) For $a \ge 2$,
the rotation number satisfies
\beql{eq307}
r(S_{ab} ) = 0 ~.
\eeq
\end{theorem}

\paragraph{Proof.}
We first consider the lower bounds in (i)--(iii).
Now $r(S_{ab} )$ is nonincreasing in $b$, and $S_{ab}$
has a fixed point whenever $b \ge 2$;
the existence of a fixed point implies $r(S_{ab}) =0$.
This establishes the lower bound and shows it is attained for all real $a$.

For the upper bound in case (iii), if $a \ge 2$ then $r(S_{ab}) =0$ 
since $r(S_{ab} ) = r(S_{ba} )$ by the conjugacy
\eqn{eq204}, and $r(S_{ba} ) =0$.

For the upper bound in the other cases, these will follow 
if we show they are attained for all sufficiently large negative $b$.
Suppose first that $-2 \le a \le 2$, and define $\theta$ so that 
$a = 2 \cos \theta$, $0 \le \theta \le \pi$.
Let $x_i(\bv)=(T_{ab}^{(i)} \bv)_x$.  Then we claim that, for $b$
sufficiently small, there is a periodic point $\bv$ of $S_{ab}$ such
that $\bv_x<0$, $x_i(\bv)> 0$, $1\le i\le n$, and
$S_{ab}^{(n+1)} (\bv)=(\bv)$, where $n$ is as in the hypothesis
($n=1$ for $a<0$).
Note that $n =\lfloor {\pi\over\theta}\rfloor.$

Now, note that a periodic point of $S_{ab}$ satisfying the constraints
must be a real eigenvector of the matrix 
$M=\left[{a\atop 1}{-1\atop 0}\right]^n\left[{b\atop 1}{-1\atop 0}\right]$.  
Since this matrix has determinant 1,
a necessary condition is that its trace be at least $2$ (if the trace were
negative, the sign condition would necessarily be violated).  This gives
a polynomial inequality in $a$ and $b$; in terms of $\theta$, this is
$$
\sin((n+1)\theta) b \ge 2 \sin(n\theta)+ 2 \sin(\theta).
$$
Note that the choice of $n$ forces $\sin((n+1)\theta)\le 0$, so we have:
$$
b\le 2 {\sin(n\theta)+\sin(\theta)\over\sin((n+1)\theta)}.
$$
Now, define
$c \ge 0$ by
$c\sin(\theta)=b \sin((n+1)\theta) -2\sin(n\theta)-2\sin(\theta)$.
 Then the eigenvector of $M$ of eigenvalue $\ge 1$ is:
$$
\bv=({\sin((n+1)\theta)\over\sin(\theta)},
         {1\over 2}(c-\sqrt{c(c+4)})+1+{\sin(n\theta)\over\sin(\theta)}).
$$
It remains only to verify that $\bv$ satisfies the sign conditions.
Now, $v_x\le 0$ by inspection, so we need consider only $x_i$, $1\le i\le n$.
By explicitly performing the appropriate matrix multiplication, we get:
$$
x_m={\sin(m\theta)\over\sin(\theta)}k+
    {\sin((n-m+1)\theta)\over\sin(\theta)},
$$
for $1\le m\le n$, where $k=1+{1\over 2}(c+\sqrt{c(c+4)})$. This is clearly
positive for all $1\le m\le n$, so the result is shown.

The remaining
upper bound case
for $a< -2$ is handled by rewriting the proof
above explicitly in terms of $a$ for the $n=1$ case;
we omit the details.~~~$\bsq$
%
%
\paragraph{Proof of Theorem \ref{th21}.}
By Theorem \ref{th32}(i) the rotation number $r(S_{ab})$ 
is nonincreasing in $b$, and by Theorem \ref{th32} 
it has rotation number attaining the endpoint values of 
the given intervals, in cases (i)--(iii).
It remains only to show that for each $a$, $r(S_{ab} )$ 
is continuous in $b$, for it must then fill out the whole interval.

We first show that for varying $b$, $r(S_{ab} )$ assumes 
every rational value in the given interval.
The proof of Theorem \ref{th32} shows that $S_{ab} (0)$ 
moves clockwise continuously in $a$ as $a$ increases.
It follows that $S_{ab}^{(n)} (0)$ has the same behavior for each $n$.
If $r = \frac{p}{q}$ is interior to the interval,
 then examining $S_{ab}^{(q)} (0)$ as $b$ increases,
its continuous variation in $b$ requires it to pass through a fixed point
$S_{ab}^{(q)} (0) =0$ such that $\theta =0$ is a 
periodic point of rotation number $\frac{p}{q}$, 
because if it does not, the rotation number never 
increases beyond $\frac{p}{q}$, a contradiction.

Any nondecreasing function $f:\RR \to [a,b]$ 
which takes every rational value in $[a,b]$ is continuous.
The nondecreasing property guarantees that the left limit
$$\lim_{x \searrow x_0^+} f(x) = f^+ (x_0)$$
and right limit
$$\lim_{x \nearrow x_0^-} f(x) = f^- (x_0)$$
both exist, with $f^- (x_0) \le f^+ (x_0)$.
If $f^- (x_0) < f^+ (x_0)$, then 
since $f(x_0)$ can take only one rational value in the interval
$(f^- (x_0), f^+ (x_0))$,
some rational value is omitted from the range of $f(x)$, a contradiction.
Thus $f^- (x_0) = f^+ (x_0)$ and $f(x)$ is continuous at $x_0$.~~~$\bsq$
%
%
\paragraph{Proof of Theorem \ref{th22}.}
Let $\mu$ be fixed.
Theorem \ref{th32}(i) shows that $r(S_{\mu \nu} )$ 
is nondecreasing in $\nu$, and $r(S_{\mu \nu}) \in [0, \frac{1}{2}]$ 
by Theorem \ref{th34}.
We have $r(S_{\mu \nu} )=0$ when $\nu \ge 2- \mu$ 
by Theorem \ref{th34}(iii) since
$a = \mu + \nu \ge 2$.
For $\nu$ negative and large the matrix 
$M = \left[ \begin{array}{cc}
b&-1 \\1 &0 \end{array}\right] 
\left[ \begin{array}{cc} a&-1 \\ 1&0 \end{array}\right]$ 
has real eigenvalues with eigenvector of sign $(+, -)$ hence
$S_{\mu \nu}$ has periodic point of period 2 and $r(S_{\mu \nu} ) =2$.

To establish continuity of $r(S_{\mu \nu} )$ in $\nu$, 
it suffices to show that $S_{\mu \nu}$ has periodic points 
with rotation number $\frac{p}{q}$ for all rational values
$0 < \frac{p}{q} < \frac{1}{2}$.
This follows similarly to the proof of Theorem \ref{th21}.
(We will later prove a stronger result in Lemma \ref{le41}.)~~~$\bsq$ \\

In the special cases where
$S_{ab}^{(n)}(0, 1) = (0, \pm 1)$ or
$S_{ab}^{(n)}(0, -1) = (0, \pm 1)$
we get significant information  about the 
behavior of $T_{ab}$.
%

\begin{theorem}~\label{th34}
Let $n\ge 1$ and suppose that
$S_{ab}^{(n)}(0,1) = (0, \pm 1),$
or $S_{ab}^{(n)}(0,-1) = (0, \pm 1),$
Then one of the following two relations holds:
\beql{eq400a}
T^{(n)}_{ab}(0,1)=(0,\lambda),
\eeq
\beql{eq400b}
T^{(n)}_{ab}(0,-1)=(0, - \lambda^{-1}).
\eeq
where $\lambda$ is a nonzero real number.

(i) If $\lambda>0$, then 
both relations above hold. In addition, 
\beql{307a}
T^{(n)}_{ab}(-1,0) =  (-\lambda,0) ~~~~\mbox{and}~~~~
T^{(n)}_{ab}(1,0) =  (\lambda^{-1},0).
\eeq
The rotation number $r(S_{ab})$ is rational.

(ii) If $\lambda < 0 $, then necessarily $\lambda=-1$.
In the first case
\beql{308a}
T^{(n)}_{ab}(0, 1)   =  (0,-1) ~~~~\mbox{and}~~~~
T^{(n)}_{ab}(- 1,0)  =  (1,0),
\eeq
while in the second case,
\beql{309a}
T^{(n)}_{ab}(0,-1)  =  (0,1) ~~~~\mbox{and}~~~~
T^{(n)}_{ab}(1,0)   =  (-1,0).
\eeq
The rotation number $r(S_{ab})$ can be irrational or rational.
\end{theorem} 

\noindent\paragraph{Proof.}
We check that the involution $R (x, y)= (y,x)$ has
\beql{eq401a}
T_{ab}^{-1} (x,y) = R^{-1} \circ T_{ab} \circ R ~,
\eeq
We have
$T_{ab}(x,y)= (ax - y, x)$ if $x \ge 0$,
and $(bx -y, x)$ if $x < 0.$ This gives 
$$
T_{ab}(x,y)=(z,w)~~ \Rightarrow~~ T_{ab}(w,z)=(y,x),
$$
where in fact $w=x$, which implies \eqn{eq401a}. 
By induction on $n \ge 1$ this yields the implication
\begin{equation}~\label{implies}
T^{(n)}_{ab}(x,y)=(z,w)~~ \Rightarrow~~  T^{(n)}_{ab}(w,z)=(y,x).
\end{equation}

Now suppose the first relation \eqn{eq400a} holds.
This yields
$$
T^{(n)}_{ab}(-1,0)= T^{(n+1)}_{ab}(0,1)
=T_{ab}(0,\lambda)
=(-\lambda,0).
$$
Using (\ref{implies}) applied to the first relation and to this 
relation yields 
\begin{eqnarray}
T^{(n)}_{ab}(\lambda,0) & = & (1,0),\\
T^{(n)}_{ab}(0,-\lambda) & = &(0,-1).
\end{eqnarray}
If $\lambda >0 $ then rescaling these by a factor $\lambda^{-1}$ gives (i),
e.g. 
$$
T^{(n)}_{ab}(0,-1)= \lambda^{-1}T^{(n)}_{ab}(0,-\lambda) 
=(0,-\lambda^{-1}).
$$
If  $\lambda <0$, then we obtain:
$$
(-\lambda,0)
=
T^{(n)}_{ab}(-1,0)
=
|\lambda|^{-1} T^{(n)}_{ab}(\lambda,0)
=
(-\lambda^{-1},0)
$$
and thus $\lambda=\lambda^{-1}=-1$. By hypothesis
$T^{(n)}_{ab}(0, 1)= (0, -1)$
and applying $T_{ab}$ to both sides gives the other relation.

Now suppose the second relation \eqn{eq400b} holds. This yields
$$
T^{(n)}_{ab}(1,0)=T^{(n+1)}_{ab}(0,-1)=T_{ab}(0,-\lambda^{-1})
=(\lambda^{-1},0).
$$
Using (\ref{implies}) applied to these relations yields
\begin{eqnarray}
T^{(n)}_{ab}( -\lambda^{-1}, 0) &= & (0,-1  )\\
T^{(n)}_{ab}(0,\lambda^{-1}) & = & (0,1).
\end{eqnarray}
If $\lambda > 0$ rescaling these equalities by $\lambda$ gives (i).
The case $\lambda < 0$ giving
the second case in (ii).
is done similarly to the first case.

If $\lambda >0$ in the first or second relation then  
$S_{ab}^{(n)}(0, 1) = (0,1)$ (resp. $S_{ab}^{(n)}(0, -1) = (0,-1)$)
so the rotation number $r(S_{ab})$ is  rational. 

If $\lambda < 0$
in these relations there is no constraint on the rationality of 
the rotation number  $r(S_{ab})$. 
In part II we give a case (ii) example having with  irrational
rotation number (\cite[Example 4.1]{LR02aa}), and also
a case (ii) example with rational rotation number
(\cite[Example 4.3]{LR02aa}). ~~~$\bsq$ \\

There are several special cases in which one can give explicit formulae for
the rotation number
$r(S_{ab})$. These are of interest in part because we have no formulae
for generic $a$ and $b$ and in part because they
 also provide several examples of
parameter values for which $T_{ab}$ is periodic.
These results all appear in Beardon, Bullett and Rippon \cite{BBR95}.
%
%
\begin{exam}
\label{ex31}
If $a=b=2\cos(\theta)$,
with $0 < \theta < \pi$,
then $r(S_{ab})={\theta\over 2\pi}$.
\end{exam}

\paragraph{Proof.}
Note that $T_{ab}$ is a linear map, since $a=b$.  If we conjugate
$T_{ab}$ by the matrix
$M=\left[{1\atop \cos(\theta)}{0\atop \sin(\theta)}\right]$, 
then we get the matrix
$\left[{\cos(\theta)\atop \sin(\theta)}
{-\sin(\theta)\atop\cos(\theta)}\right],$
which is rotation by $\theta$.
 Since $M$ is orientation-preserving for $0 < \theta < \pi$, conjugation
by $M$ preserves the rotation number of $S_{ab}$; rotation by $\theta$ has
rotation number ${\theta\over 2\pi}$, thus the formula given.  Note that
when ${\theta\over 2\pi}$ is rational, $T_{0\nu}$ is periodic,
and the period is computed from the rotation number.~~~$\bsq$
%
%
\begin{exam}\label{ex32}
If $a=2\cos({\pi\over n})$, for $n \ge 2$, and $b=2\cos( \theta)$ 
with $0 < \theta < \pi$, then
$r(S_{ab})={\theta\over \pi+n\theta}$. 
When $\frac{\theta}{2\pi} = \frac{p}{q}$
is rational, then $T_{ab}$ is periodic with period $2np + q$, for $q$ odd,
and half this if $q$ is even.
\end{exam}

\paragraph{Proof.}
Here, we have a slightly more complicated situation, in that
the map $T_{ab}$  is not linear.  Again, we conjugate by
$M=\left[{1\atop \cos(\theta)}{0\atop \sin(\theta)}\right]$,
to get 
$$
T'(\bv)=
\cases{
\bv\left[{\cos(\theta)\atop \sin(\theta)}
{-\sin(\theta)\atop\cos(\theta)}\right]
 ,&if $(\cos(\theta),\sin(\theta))\cdot\bv\le 0$,\cr
\bv M\left[{2\cos({\pi\over n})\atop 1}{-1\atop 0}\right] M^{-1},&
if $(\cos(\theta),\sin(\theta))\cdot\bv\le 0$.\cr}
$$
The important thing to note is that the second case will always occur
exactly $n$ times in a row, and that the $n$th power of that matrix is
$-1$.  If, in formula \eqn{eq201A}, we consider only the subsequence of
iterates not in the middle of such a block, we get the following:
$$
r(S_{ab})= \lim_{i\to\infty} {i{\theta\over 2\pi}+
                {1\over 2}\lfloor i{\theta\over\pi}\rfloor\over
                   i+\lfloor i{\theta\over\pi}\rfloor n}.
$$
In the limit, we can ignore the floors, to get ${\theta\over\pi+n\theta}$,
as claimed.  Here, again, when ${\theta\over 2\pi}$ is rational, $T_{ab}$
is periodic.~~~$\bsq$ \\

As a special case of Example \ref{ex32}, 
in $(\mu, \nu)$ parameters,
taking   $\nu=0$,   $\mu = 2\cos \frac{\pi}{n}$,
for $n \ge 3$, yields a $T_{\mu\nu}$ which
is periodic with period $n^2.$ This corresponds to the parameters
$a=2\cos \frac{\pi}{n}$, 
$b = 2 \cos \frac{(n-1)\pi}{n}$ in Example \ref{ex32}.

%
%
\begin{exam}\label{ex33}
For each integer $n \ge 2$,
if  both $a,b\le 0$, and $ab=4\cos^2({\pi\over 2n})$,
then $r(S_{ab})={2n-1\over 4n}$, and $T_{ab}$ is periodic.
\end{exam}

\paragraph{Proof.}
For this case, the simplest approach is to explicitly calculate
$T_{ab}^{(m)} (\bv)$, for $1\le m \le 4n$, and plug that information into
equation \eqn{eq201A}.  It is most convenient to take $\bv_x,\bv_y <0$.  
We can now calculate:
\begin{eqnarray*}
x_{2i}&=&
{
  \sin({(i+1)\pi\over n})+\sin({i\pi\over n})
\over\sin({\pi\over n})}\bv_x-
{\sin({i\pi\over n})\over\sin({\pi\over n})}a \bv_y\le 0,\\
x_{2i+1}&=&
{
  \sin({(i+2)\pi\over n})+2\sin({(i+1)\pi\over n})+\sin({i\pi\over n})
\over
  \sin({\pi\over n})
} {\bv_x\over a}-
{\sin({(i+1)\pi\over n})+\sin({i\pi\over n})
 \over\sin({\pi\over n})}\bv_y\ge 0,\cr
\end{eqnarray*}
for $0\le i\le (n-1)$.  One can now calculate that
$T_{ab}^{(2n)} (\bv) =-\bv$; since, by symmetry, we could as
easily have started with $\bv_x,\bv_y>0$, we can conclude that
$T_{ab}^{(4n)}=1$.  By inspecting the sequence of quadrants we pass
through in a period, we conclude that in $4n$ steps, we wrap around the
origin $2n-1$ times; $r(S_{ab})$ is thus $\frac{2n-1}{4n}$.~~~$\bsq$

Example \ref{ex33} gives a countable collection of one-parameter
families, on
which $T_{ab}$ is a periodic map ; the question of whether other any other
continuous families of periodic parameter values exist is still open. 
Example 4.1 in part II gives a countable number of periodic maps
with isolated periodic parameter values. It gives a 
one-parameter family with fixed symbolic dynamics, where the
rotation number changes in the family. The fixing of the symbolic
dynamics bt $T^{(8)}_{ab}(0,-1)=
(0, 1)$ yields a one-parameter family 
and since the rotation number varies in the 
family,  the rational rotation number
points, which give periodic maps, must have
isolated parameter values. 
%
%
%
%
%

\section{Rational Rotation Number}
\hsp
We first characterize those cases when $T_{ab}$ is a periodic map.

\paragraph{Proof of Theorem \ref{th23}.}
If $T_{ab}$ is periodic, then $(0,1)$ is a periodic point.

Suppose, conversely, that $(0,1)$ is a periodic point of $T_{ab}$ of minimal
period $p$.
Then so is $(1,0)$, with the same period $p$, because $T_{ab}^{-1}$ 
is obtained from $T_{ab}$ by interchanging $x$ and $y$ coordinates,
cf. Theorem~\ref{th34}.
Next $T_{ab} ((0,-1)) = (1,0)$ hence $(0,-1)$ is a periodic point.
We use the piecewise linear nature of $T_{ab}$ to prove that $T_{ab}$ 
is periodic.
An iterate  $T_{ab}^{(n)}$ is linear on a sector $\Sigma$ as 
long as the $x$-axis never intersects the interior of 
$T_{ab}^{(j)} (\Sigma )$ for $0 \le j < n$.
Because $(0,1)$ and $(0,-1)$ are both periodic, the image of the $x$-axis
under all iterates $\{T^{(n)} : n > 0 \}$ splits the plane into a finite
number of sectors, and the image of each of these 
sectors under $T$ is another sector.
It follows that $T_{ab}^{(p)}$ permutes these sectors 
and is linear on each one of them.
But each sector contains two linearly independent fixed points 
of $T_{ab}^{(p)}$, from which it follows that $T_{ab}^{(p)}$ 
is the identity map on each sector.
Since these sectors partition the entire plane, $T_{ab}$ 
is periodic with period $p$.~~~$\bsq$ \\

The parameter values giving periodic maps $T_{ab}$ 
can be classified by the period length $p$, and by the symbolic dynamics
of the iterates, in which the
symbolic dynamics describes the $x$-coordinate signs of each iterate
as belonging to  $\{+, - , 0\}$.
There are  a countable collection of such data. In each case,
either one gets (one or two) isolated parameter values, or else a 
one-parameter family of values $(a,b)$; the family is 
determined by the requirement that
$T_{ab}^{(p)}(0,1) = (0,1),$ cf. Theorem~\ref{th34}.
Example~\ref{ex33} exhibited a countable number of one-parameter
families of periodic maps.
There are a countable number of
 parameter values that are isolated 
in the sense that they are the unique point in some open
neighborhood in the parameter space for a periodic map with
the given period and symbolic dynamics. Examples appear in
part II, as described at the end of \S3. \\

  We also note that the  parameter values for periodic points can be 
subdivided into two types.

{\em Type (i)}.
$(0,1)$ and $(0, -1)$ are in the same periodic orbit of $T_{ab}$, 
necessarily of even period.

{\em Type (ii)}.
$(0,1)$ and $(0, -1)$ are in separate periodic orbits of $T_{ab}$, 
having the same (even or odd) period.

\noindent
Examples show that both  types of periodic orbit occur. \\

To prove Theorem \ref{th24}, which
characterizes the possible dynamics of $T_{ab}$
when the rotation number is rational, 
 we use the $(\mu,\nu )$-parameter space.
The following auxiliary lemma is a main tool in proving Theorem~\ref{th24};
its important feature is the uniqueness assertion.

%
%
\begin{lemma}\label{le41}
For every rational $r= \frac{p}{q}$ with $0 < r < \frac{1}{2}$ and each fixed
$\mu \in \RR$, the following holds.
For each $\bv = e^{i \theta} \in S^1$ there is a unique parameter value
$\nu \in \RR$ such that $\bv$ is a periodic point of $S_{\mu \nu}$ 
of rotation number $r$.
\end{lemma}

\paragraph{Proof.}
Let $\mu$ be fixed.
We first show that for fixed $\theta$, the second iterate 
$S_{\mu \nu}^{(2)} ( \theta )$ moves strictly clockwise as $\nu$ increases.
The proof of Theorem \ref{th32}(ii) shows that 
$\frac{d}{d\nu} S_{\mu \nu}^{(2)} (\theta ) \le 0$.
If $x_0 = \bv_x$ and $x_1 = (S_{\mu \nu} (\bv ))_x$ then for 
$\frac{d}{d\nu} S_{\mu \nu}^{(2)} (\bv ) =0$ it is necessary 
that $x_0 = x_1 =0$.
But if $x_0 = 0$, then $x_1 = - \bv_y$ which is impossible for 
$\bv \in S^1$.
By the same argument, for each $n \ge 2$, and fixed $\theta$, 
$S_{\mu \nu}^{(n)} (\theta)$ moves strictly clockwise as $\nu$ increases.

The families of Example \ref{ex33} give
us some $\nu_0$ such that $\bv$ is
a periodic point of $S_{\mu\nu_0}$ of rotation number ${2q-1\over 4q}$;
consider, now, the behavior of $S^{(4q)}_{\mu\nu} (\bv )$ as $\nu$ 
changes from $\nu_0$.
If $\nu$ decreases from $\nu_0$, $r(S_{\mu\nu})>{2q-1\over 4q}>r$, 
so we need
consider only $\nu>\nu_0$.  As we increase $\nu$, 
$\bw_{\nu} := S^{(4q)}_{\mu\nu} (\bv )$ moves
smoothly (and strictly) clockwise.
If we view
this as a function of $\nu$ to the point $\theta_\nu$ with
$\bw_\nu = e^{2 \pi i \theta_{\nu}}$, we can 
 take the lift 
$\tilde S (\nu) = \theta_\nu $ of
this function to ${\bf R}$, with the initial condition 
$\tilde S (\nu_0)= 2q-1$; this is
continuous, strictly decreasing, and for large enough $\nu$ 
attains a value below $4$.
($\tilde{S} (\nu)\ge 4$ 
implies $r(S_{\mu\nu})\ge {1\over q}$.)
Since $2q-1 > 4p \ge 4$, there is a
unique $\nu$ such that $\tilde S(\nu)=4p$.  It is then easy to see that
$\bv$ is a periodic point of $S_{\mu\nu}$ of rotation number
$r={p\over q}$.~~~$\bsq$

%
%

\paragraph{Proof of Theorem \ref{th24}.}
We use the $(\mu,\nu)$-parameter space.  We fix $\mu$, and allow $\nu$ to
vary; we can use the monotonicity in $\nu$ to deduce the behavior for every
$\nu$ (for that rotation number) from the behavior from one special $\nu$.
In particular, let $\nu_0$ be the unique $\nu$ such that $(0,1)$ is a
periodic point of $S_{\mu\nu_0}$ of rotation number $r(S_{\mu\nu})={p\over
q}$.  Now, if $(0,1)$ is a periodic point of $T_{\mu\nu_0}$, then
$T_{\mu\nu_0}$ is periodic, by Theorem \ref{th23}; since rational rotation
number implies existence of a periodic point, the uniqueness part of Lemma
\ref{le41} then tells us that $\nu=\nu_0$ is the only value of $\nu$ for
which $r(S_{\mu\nu})=\frac{p}{q}$. This situation  corresponds to case
(iii) of the theorem, because  we will show below that in all remaining
cases $S_{ab}$ has at most two periodic orbits.

Thus, in what follows we  assume that $(0,1)$ is not a periodic point of
$T_{\mu\nu_0}$, so the map
is not a periodic map. 
  Then there must be some positive real $\lambda \ne 1$ 
such that
$T^{(q)}_{\mu\nu_0}(0,1)=(0,\lambda)$. 
We can, without loss of generality, assume that $\lambda>1$; 
the proof for $\lambda<1$ is analogous. For
the orbit of $(0, -1)$, we also have by Theorem~\ref{th34}(i)
that $T^{(q)}_{\mu\nu_0}(0,-1)=(0,-\lambda^{-1})$.

     Note first that if we plot the orbits of $(0,1)$ and $(0,-1)$ under
$S_{\mu\nu_0}$, then the two orbits are alternating as we move around the
circle.  Also, the sector bounded by a point of one orbit and a neighboring
point from the other orbit is linear under the action of 
$T^{(q)}_{\mu\nu_0}$.
Now in each sector $T^{(q)}_{\mu\nu_0}$ has two eigenvectors, given
by the endpoints of the sector, so it can have no other eigenvector
in the sector, unless it is a multiple of the identity
on the whole sector. But this cannot happen, for if $T^{(q)}_{\mu\nu_0}$ 
were a multiple of the identity on the sector, necessarily $\pm I$,
then  each endray would consists of periodic points of period $q$ or $2q$, 
and the orbit of some point on 
one of these endrays contains $(0,1)$, so $(0,1)$ would be a periodic point,
contradicting our assumption. We conclude that 
$S_{\mu\nu_0}$ has no other periodic orbits.

Now, for any $0\le n<q$, let $\bv_0(\nu)\in S^1$ be
$S^{(-n)}_{\mu\nu}(0,1)$, and let $\bv_1(\nu)$ and $\bv_{-1}(\nu)$
be the points of the form $S^{(-i)}_{\mu\nu}(0,-1)$, $0\le i<q$ immediately
counterclockwise and clockwise of $\bv_0(\nu)$, respectively.  For
every point in the interval from $\bv_{-1}(\nu)$ to $\bv_1(\nu)$, 
label the point
as {\it clockwise} if $S^{(q)}_{\mu\nu}$ 
moves the point clockwise, and similarly
label it as {\em counterclockwise} if $S^{(q)}_{\mu\nu}$ 
moves the point
counterclockwise (fixed points will receive no label).  Now, at
$\nu=\nu_0$, the points between $\bv_{-1}$ and $\bv_0$ will be 
labelled clockwise,
and the points between $\bv_0$ and $\bv_1$ will be labelled counterclockwise,
since $\bv_0$ is an attracting fixed point of $S^{(q)}_{\mu\nu}$, and
$\bv_{\pm 1}$ are
repelling fixed points of $S^{(q)}_{\mu\nu}$.  As $\nu$ decreases,
 every point labelled
counterclockwise will remain labelled counterclockwise (by the proof of
Theorem \ref{th22}); further, $\bv_0$ and $\bv_{\pm 1}$ will be labelled
counterclockwise.  Now, $T^{(q)}_{\mu\nu}$ is linear on the interval from
$\bv_0$ to $\bv_1$, so $S^{(q)}_{\mu\nu}$ can have at most 
two fixed points in
the interval.  As long as there are points labelled clockwise in the
interval, $S^{(q)}_{\mu\nu}$ must have exactly two fixed points; these must
move continuously toward each other; when the two fixed points coincide, we
are left with one fixed point; after that point, there can be no fixed points
in the interval.  Thus, we can conclude that the interval
from $\bv_{-1}$ to $\bv_1$ always has either two, one, or zero 
fixed points of
$S^{(q)}_{\mu\nu}$, in the range of $\nu\le\nu_0$ such that
$r(S_{\mu\nu})={p\over q}$; the proof for $\nu\ge\nu_0$ is analogous.
Now, the proof applies for any $\bv_0$; there are exactly $q$ 
possible choices
of $\bv_0$, and the intervals from $\bv_{-1}$ to $\bv_1$ around 
these choices of $\bv_0$ exhaust $S^1$; therefore, we can conclude 
that the number of distinct periodic
orbits of $S_{\mu\nu}$ must equal the number of distinct fixed points of
$S^{(q)}_{\mu\nu}$ between $\bv_{-1}$ and $\bv_1$, which is at most two.
There are then three cases: \\

{\em Case 1}: {\em $S_{\mu\nu}$ has no periodic orbit}. \\

Then $r(S_{\mu\nu})$ is irrational, contradicting the hypothesis
that it is rational; this case cannot occur. \\

{\em Case 2}: {\em $S_{\mu\nu}$ has exactly two periodic orbits}. \\ 

Let $\bw$ be a fixed point
corresponding to one of the orbits, and $\bw'$ a fixed point corresponding
to the other orbit.  Then there must be some $\lambda$ such that
$T^{(q)}_{\mu\nu}(\bw)=\lambda \bw$, 
and $T^{(q)}_{\mu\nu}(\bw')=\lambda^{-1} \bw'$.
For every point $\bv$, $S^{(i)}_{\mu\nu}(\bv)$ comes arbitrarily 
close to $\bw$
(for $i>0$); there is an interval around $\bw$ which behaves linearly under
$T^{(q)}_{\mu\nu}$, so every point eventually is transformed by the same
map.  $T^{(q)}_{\mu\nu}$ has an eigenvector of eigenvalue $>1$; it is easy
to see that every point must diverge exponentially except for the
eigenvector of eigenvalue $<1$, which converges exponentially to $0$.
Similarly, in the other direction, every point diverges exponentially except
for the eigenvector of eigenvalue $>1$, which 
converges exponentially to $0$. This is case (ii) in the theorem. \\

{\em Case 3}: {\em $S_{\mu\nu}$ has exactly one periodic orbit}. \\

If we reflect a periodic
orbit of $S_{\mu\nu}$ around the line $x=y$, we must get a periodic orbit;
thus, in this case, the periodic orbit must be symmetrical around $x=y$.
It follows, then, that the orbit lifts to a periodic orbit of $T_{\mu\nu}$
($T^{(q)}_{\mu\nu}$ and $T^{(-q)}_{\mu\nu}$ must multiply the vectors in the
orbit by the same amount).  In a sufficiently small neighborhood of a point
of the orbit, $T^{(q)}_{\mu\nu}$ is linear; it has exactly one eigenvector
of eigenvalue one, but has eigenvalue one with multiplicity two; it follows
that every point in that neighborhood diverges linearly.

     Every point on the circle has a value of $\nu$ associated to it by
the above reasoning such that it is a periodic point of $S_{\mu\nu}$.
By Lemma \ref{le41}, there can be no other $\nu$ 
such that $S_{\mu\nu}$ has a periodic
point of rotation number $\frac{p}{q}$.  This
establishes case (i) of the theorem. ~~~$\bsq$ \\

%
%
%
%
%

\section{Irrational Rotation Number}
\hsp

When the rotation number is irrational, the analysis
of the dynamics  becomes
significantly more complicated. 
On the level of the circle maps
$S_{ab}$, their smoothness properties
lead to the following result.
\setcounter{figure}{0}
%
%

\begin{theorem}
\label{th51}
If the rotation number $r(S_{ab})=r$ is irrational,
 then there exists
a homeomorphism $h:S^1\to S^1$ such that 
$h \circ S_{\mu\nu} \circ h^{-1}=\Theta_r$,
where $\Theta_r:S^1\to S^1$ is rotation by $r$
($\theta\mapsto\theta+2\pi r$).
\end{theorem}

\paragraph{Proof.}
By Theorem \ref{th31} ${d\over d\theta} S_{ab}(\theta)$ is
continuous and of bounded variation.
The topological conjugacy result then follows from Denjoy's
Theorem, see Nitecki \cite[p. 41 and Theorem, p. 45]{Ni71},
or Herman \cite[VI.4]{He79}.~~~$\bsq$ \\

This effectively characterizes the 
behavior of $S_{ab}$, but it still leaves several possibilities for
the dynamical behavior
of $T_{ab}$.  In particular, it does not rule out the possibility 
that $T_{ab}$ has 
a weakly localized orbit (one that converges to 0 as $n \to \pm \infty$) 
or of a divergent orbit (one that becomes unbounded as $n \to \pm \infty)$.
We do not know if either of these possibilities occurs.

One major problem is that of establishing
conditions under which
$T_{ab}$ has an invariant circle.
When this holds, scale invariance implies there
is an invariant circle through every point (except the origin), 
and the dynamics is
bounded. Herman  \cite[VIII.2.4]{He79} shows that this occurs
if and only if $T_{ab}$ is topologically conjugate to a rotation
of the plane, and that this holds if and only if the map $S_{ab}$
is $C^{1}$-conjugate to a rotation of the circle. Herman  conjectures
that there do exist maps $T_{ab}$ 
with $S_{ab}$ having irrational rotation such that $T_{ab}$
is not topologically conjugate to a rotation of the plane.

We give the following condition for the existence of an invariant
circle. 
%
%
\begin{theorem}\label{th52}
Let $T: \RR^2 \to \RR^2$ be any homeomorphism which is scale-invariant;
that is
$T(\lambda \bv ) = \la T ( \bv )$ for all $\lambda \ge 0$, 
all $\bv \in \RR^2$.
Let $S(\bv ) = \frac{T(\bv )}{\|T( \bv ) \|}$ be its 
associated circle map, and suppose that $S$ is 
topologically conjugate to an irrational rotation.
If $T$ has a forward orbit $\sO^+ (\bv_0 ) = \{\bv_n : n \ge 0\}$ such that
$$\frac{1}{C} \le \| \bv_n \| \le C , \qquad \mbox{all} \quad n \ge 0~,$$
for some positive constant $C>1$, then $T$ has an invariant circle, 
and all orbits are bounded.
\end{theorem}

\paragraph{Proof.}
For each $n\ge 0$, let $L_n$ be the line segment from $0$ to $T^{(n)}(\bv)$,
and let
$$
I_m=\bigcup_{n\ge m} L_n, \quad\mbox{and}\quad
I'=\bigcap_{m\ge 0} \overline{I_m}.
$$
Now, $T(I_m)=I_{m+1}$, so $T(\overline{I_m})=\overline{I_{m+1}}$, 
and $\overline{I_m}$ is a
decreasing sequence of sets.  Therefore, we can conclude that $I'$ is an
invariant set of $T$.  Furthermore, $\overline{I_m}$ 
contains a ball of radius $1/C$ around
$0$ for each $m$, since $S^{(n)} (\bv)$ is dense in $S^1$, so $I'$ contains
a ball of radius $1/C$.

     Now, consider the set $B:=\partial I'$, the boundary of $I'$.  
$B$ is also
an invariant set of $T$ which does not contain $0$, and it contains 
at least one
point in every direction from $0$.  If we can show that $B$ is homeomorphic
to $S^1$, we will be done.  It suffices to show that $B$ is connected, and
contains at most one point in every direction from $0$; projection onto the
unit circle then provides a homeomorphism to $S^1$.  Connectedness of $B$
follows from the fact that $I'$ is closed, bounded, and star-convex, so
we have left only to show that every ray from the origin hits $B$ in exactly
one point.

Suppose that $\bv_1,\bv_2\in B$, with $\bv_2=k \bv_1$, $k>1$.  
Then, for all $n\ge 0$,
$$
{k\over C}\le |T^{(n)}(\bv_2)|\le C \,,
$$
and $|T^{(n)}(\bv_2)|=k|T^{(n)}(\bv_1)|$.  Now, $B$ contains the 
line segment 
from $\bv_1$ to
$\bv_2$, and thus contains the line segment from $T^{(n)}(\bv_1)$ to 
$T^{(n)}(\bv_2)$ for
all $n\ge 0$.  These line segments have length bounded from below,
all lie within the closed ball of radius $C$ around $0$, and lie on a dense
set of rays.  It thus follows that $\bar B$ has non-empty interior.  But, $B$
is a boundary, so this is impossible.  Therefore, every ray must hit $B$ in
at most one point, so $B$ is an invariant circle.~~~$\bsq$ \\

Theorem \ref{th52} applies to $T_{ab}$, and 
gives a criterion 
when it has  invariant circles.
We conclude by establishing  the reflection
 symmetry of invariant
circles about the line $x=y$ (when they exist); 
this symmetry is manifest
in plots of invariant circles given in part II.

\begin{theorem}~\label{th52a}
If the map $T_{ab}$ has associated
irrational rotation number $r(S_{ab})$,
they any invariant circle of  $T_{ab}$
is symmetric under the
involution $R(x,y) = (y,x).$
\end{theorem}
\noindent\paragraph{Proof.}
Let $\sC$ denote an invariant circle.
The  conjugacy \eqn{eq204a} by $R$ shows that the set
$R(\sC)= R^{-1}(\sC)$ is invariant under 
$T_{ab}^{-1}=  R^{-1} \circ T_{ab}\circ R$
and hence under $T_{ab}$, so that it too is an invariant circle
of   $T_{ab}$.
The invariant circle $\sC$ intersects the line
$x=y$, and points on this line are left fixed by $R$, so
that $\sC$ and $R(\sC)$ have at least one common point. 
For irrational rotation number, any
two invariant circles having a common point must be identical,
because they have a dense orbit in common; thus
$\sC= R(\sC)$.
~~~$\bsq$


\end{document}